# A counterexample to the "hot spots" conjecture

By Krzysztof Burdzy and Wendelin Werner*

**Abstract**

We construct a counterexample to the "hot spots" conjecture; there exists a bounded connected planar domain (with two holes) such that the second eigenvalue of the Laplacian in that domain with Neumann boundary conditions is simple and such that the corresponding eigenfunction attains its strict maximum at an interior point of that domain.

## 1. Introduction

The "hot spots" conjecture says that the eigenfunction corresponding to the second eigenvalue of the Laplacian with Neumann boundary conditions attains its maximum and minimum on the boundary of the domain. The conjecture was proposed by J. Rauch at a conference in 1974. Our paper presents a counterexample to this conjecture.

Suppose that $D$ is an open connected bounded subset of $\mathbf{R}^d$, $d \geq 1$. Let $\{\varphi_1, \varphi_2, \dots\}$ be a complete set of $L^2$-orthonormal eigenfunctions for the Laplacian in $D$ with Neumann boundary conditions, corresponding to eigenvalues $0 = \mu_1 < \mu_2 \leq \mu_3 \leq \mu_4 \leq \dots$. The first eigenfunction $\varphi_1$ is constant.

THEOREM 1. *There exists a planar domain $D$ with two holes (i.e., conformally equivalent to a disc with two slits) such that the second eigenvalue $\mu_2$ is simple (i.e., there is only one eigenfunction $\varphi_2$ corresponding to $\mu_2$, up to a multiplicative constant) and such that the eigenfunction $\varphi_2$ attains its strict maximum at an interior point of $D$.*

There remains a problem of proving the conjecture under additional assumptions on the geometry of the domain. Since our method does not seem to be able to generate a counterexample with fewer than two holes, it is natural to ask if this failure has causes of fundamental nature.

---

*Research of the first author partially supported by NSF grant DMS-9700721.

*Problem* 1. Does the "hot spots" conjecture hold in all planar domains which have at most one hole?

Using a probabilistic coupling argument, Bañuelos and Burdzy (1999) proved the conjecture for some "long and thin" (not necessarily convex) planar domains and for some convex planar domains with a line of symmetry. We know of only one other published result on the conjecture; it is contained in a book by Kawohl (1985).

We are grateful to David Jerison and Nikolai Nadirashvili for telling us about their forthcoming results. They include a proof of the "hot spots" conjecture for convex planar domains and a different counterexample.

See the introduction to Bañuelos and Burdzy (1999) for a detailed review of various aspects of the "hot spots" conjecture, and a complete reference list. Our techniques are very close to those introduced in that paper so we will be rather brief and we ask the reader to consult that paper for more details.

We would like to thank Rodrigo Bañuelos and David Jerison for very useful advice, and the anonymous referee for suggesting a short proof of Lemma 1. The second author had the pleasure of being introduced to the problem by Jeff Rauch, the proposer, at E.N.S. Paris in 1995.

## 2. Domain construction

Before describing precisely our domain $D$, let us now give a short intuitive argument that provides some heuristic insight into our counterexample. Consider a planar domain that looks like a bicycle wheel with a hub, at least three very very thin spokes and a tire. Consider the heat equation in that domain with Neumann boundary conditions and an initial temperature such that the hub is "hot" and the tire is "cold." Due to the fact that the cold arrives in the hub only via the spokes, the "hottest spot" of the wheel will be pushed towards the center of the hub. This implies that the second Neumann eigenfunction in the domain attains its maximum near the center of the hub and therefore not on the boundary of the domain.

For technical reasons that will become apparent in the proof, our domain $D$ does not quite look like a bicycle wheel, but it does have a "hub," three "spokes" and a "tire."

We will use $\underline{0}$ as an abbreviation for $(0,0)$. Let $G$ be the group (containing six elements) generated by the symmetry $s$ with respect to the horizontal axis and the rotation around $\underline{0}$ by the angle $2\pi/3$. We will use the point-to-set mapping $\mathcal{T}x = \{\sigma(x),\ \sigma \in G\}$. Typically, $\mathcal{T}x$ contains six points. The meaning of $\mathcal{T}K$ for a set $K$ is self-evident.

Suppose $\varepsilon \in (0, 1/200)$ is a very small constant whose value will be chosen later in the proof. Let us name a few points in the plane,

$$A_1 = \underline{0}, \quad A_2 = (1/7, \sqrt{3}/7), \quad A_3 = (5, 1/100),$$
$$A_4 = (11/2, 1/200), \quad A_5 = (6, \varepsilon), \quad A_6 = (13/2, 1/200),$$
$$A_7 = (7, 1/100), \quad A_8 = (8, 8\sqrt{3}), \quad A_9 = (9, 9\sqrt{3}),$$
$$A_{10} = (235, 0).$$

Let $D_1$ be the domain whose boundary is a polygon with consecutive vertices $A_1, A_2, A_3, A_4, A_5, A_6, A_7, A_8, A_9, A_{10}$ and $A_1$. Let $D_2$ be the closure of $\mathcal{T}D_1$ and let $D_3$ be the interior of $D_2$. Finally, we obtain $D$ be removing the line segment between $(-18, 0)$ and $(-16, 0)$ from $D_3$. We will show that $D$ has the properties stated in Theorem 1.

The domain $D_3$ has three holes while $D$ has only two, because of the cut between $(-18, 0)$ and $(-16, 0)$.

Let $\alpha_1$ and $\alpha_2$ be the minimum and maximum of the angles between vectors $\overrightarrow{A_j A_{j+1}}$, $j = 1, 2, \ldots, 9$, and the horizontal axis. We have chosen the points $A_j, j = 1, \ldots, 10$, in such a way that $\alpha_2 - \alpha_1 < \pi/2$; this fact will be useful at the end of the proof, when we apply results of Bañuelos and Burdzy (1999).

## 3. The second eigenvalue is small

In this section, we will prove the following result.

LEMMA 1. *For every $\delta > 0$, there exists $\varepsilon_0 > 0$ such that for all $\varepsilon \in (0, \varepsilon_0)$, the second Neumann eigenvalue $\mu_2$ in the domain $D = D(\varepsilon)$ defined in Section 2 is not greater than $\delta$.*

*Proof.* Recall the points $A_j$ defined in Section 2. Let $A_{11} = (6 + 100\varepsilon/(1 - 200\varepsilon), 0)$; this point lies at the intersection of the line containing $A_4$ and $A_5$ and the horizontal axis. Suppose that $\varepsilon < 1/1600$. Let $d_1(z)$ denote the Euclidean distance between $z$ and $A_{11}$ and define for all $z \in D_1$, the function $f_1$ as follows.

- $f_1(z) = 0$ if $|z| > 6$ or if $d_1(z) < 400\varepsilon$,
- $f_1(z) = \log(d_1(z)/400\varepsilon)/\log(1/800\varepsilon)$ if $|z| \leq 6$ and $d_1(z) \in [400\varepsilon, 1/2]$,
- $f_1(z) = 1$ if $|z| < 6$ and $d_1(z) > 1/2$.

Extend $f_1$ into a continous function on $D$ invariant under $G$. The function $f_1$ is equal to 1 in the inside region of $D$, it is equal to 0 in the exterior region and it slopes from 1 to 0 in the inside parts of the narrow channels connecting the two regions. Note that $f_1$ satisfies the Neumann boundary conditions on $\partial D$.

It is clear that when $\varepsilon \to 0$, the two quantities $\int_D |f_1|$ and $\int_D |f_1|^2$ remain bounded and bounded away from 0, and it is elementary to check that $\int_D |\nabla f_1|^2 \to 0$ when $\varepsilon \to 0$.

Similarly, define a function $f_2$ that is equal to 1 in the exterior domain, to 0 in the inside part of $D$ and that slopes to 0 in the exterior parts of the narrow necks connecting the two. More precisely, let $A_{12}$ lie at the intersection of the horizontal line and the line containing $A_5$ and $A_6$, let $d_2(z) = \text{dist}(z, A_{12})$ and for $z \in D_1$,

- $f_2(z) = 0$ if $z < 6$ or if $d_2(z) < 400\varepsilon$,
- $f_2(z) = \log(d_2(z)/400\varepsilon)/\log(1/800\varepsilon)$ if $|z| > 6$ and $d_2(z) \in [400\varepsilon, 1/2]$,
- $f_2(z) = 1$ if $|z| > 6$ and $d_2(z) > 1/2$.

Extend $f_2$ in a continuous fashion to $D$ so that it is invariant under the action of $G$. Both $\int_D |f_2|$ and $\int_D |f_2|^2$ remain bounded and bounded away from 0 when $\varepsilon \to 0$. It is easy to check that $f_2$ satisfies the Neumann boundary conditions and that $\int_D |\nabla f_2|^2 \to 0$ when $\varepsilon \to 0$.

Finally, define for all $\varepsilon > 0$, the function $f$ on $D$ by

$$f(z) = \frac{f_1(z)}{\int_D f_1} - \frac{f_2(z)}{\int_D f_2}.$$

The function $f$ satisfies the Neumann boundary conditions. As the supports of $f_1$ and $f_2$ are disjoint, it is clear that $\int_D |f|^2$ remains bounded and bounded away from 0 when $\varepsilon \to 0$, and that $\int_D |\nabla f|^2 \to 0$ when $\varepsilon \to 0$. Since $f$ is orthogonal to the constant function 1 (i.e., to the lowest eigenfunction) because $\int_D f = 0$, we conclude that the second Neumann eigenvalue $\mu_2$ in $D$ satisfies

$$0 < \mu_2 \leq \frac{\int_D |\nabla f|^2}{\int_D |f|^2}.$$

Hence, $\mu_2 = \mu_2(\varepsilon) \to 0$ when $\varepsilon \to 0$.  □

## 4. Nodal line of the second eigenfunction

In this section we will show that the nodal lines of any second eigenfunction (i.e., any eigenfunction corresponding to $\mu_2$) are confined to a small subset of $D$ when $\varepsilon$ is small.

More precisely, consider any second Neumann eigenfunction $\varphi_2$ in $D$. Let $\Gamma$ be its nodal line, i.e., $\Gamma = \{x \in D, \varphi_2(x) = 0\}$ (note that the line $\Gamma$ is not necessarily connected).

Recall that $s$ denotes the symmetry with respect to the horizontal axis, and define for $j = 3, 4, 5, 6, 7$ the line segments

$$K_j = \overline{A_j s(A_j)}.$$

Let $M^o$ denote the part of $D$ between $K_3$ and $K_7$, and define $M = \mathcal{T} M^o$ (i.e. "the union of the three bridges").

The goal of this section is to prove that

LEMMA 2. *For all small enough $\varepsilon > 0$, $\Gamma \subset M$.*

We divide the proof into several steps.

*Step* 1. Let $\Gamma_1$ denote a connected component of the nodal line. Suppose that $\Gamma_1$ intersects $D \setminus M$ and that the diameter of $\Gamma_1$ is less than $10^{-10}$. We will show that this assumption leads to a contradiction, if $\varepsilon$ is sufficiently small.

As the diameter of $\Gamma_1$ is less than $10^{-10}$ and $\Gamma_1 \not\subset M$, it is easy to see that $\Gamma_1$ has to cut off a domain $D_4$ from $D$ of diameter less than $10^{-6}$ (the boundary of $D_4$ would consist of $\Gamma_1$ and a piece of $\partial D$). It is also easy to prove that the first eigenvalue $\lambda_1$ for the mixed problem in $D_4$, with the Dirichlet conditions on $\Gamma_1$ and the Neumann conditions elsewhere on $\partial D_4$, is larger than some $\lambda_0 > 0$, independent of $\varepsilon < 1/200$ and the shape and location of $\Gamma_1$ (but using the fact that $\Gamma_1$ intersects $D \setminus M$). Since $\lambda_1 = \mu_2$, we can adjust $\varepsilon$ to make $\mu_2 < \lambda_0$ using Lemma 1, and we can thus rule out the possibility that the diameter of $\Gamma_1$ is less than $10^{-10}$.

*Step* 2. We now collect some simple facts on reflected Brownian motion in $D$. We define some further sets: $D \setminus M$ consists of two connected components, the inner one $I$ (the one containing $\underline{0}$) and the exterior one $E$. Also, let $M_i^o$ ($M_e^o$) denote the part of $D$ between the line segments $K_3$ and $K_5$ ($K_5$ and $K_7$). Put $M_e = \mathcal{T} M_e^o$ and $M_i = \mathcal{T} M_i^o$.

In the rest of the paper, $X_t = (X_t^1, X_t^2)$ will denote reflected Brownian motion in $D$ (with normal reflection on $\partial D$). For all $U \subset D$, $\tau_U$ will denote the first hitting time of $U$ by $X$; i.e.,

$$\tau_U = \inf\{t > 0 : X_t \in U\}.$$

Define $Z_t = |X_t^1 - 6|$. As long as $X_t$ stays in $M^o$, the process $Z_t$ is a one-dimensional Brownian motion reflected at 0, with some local time push always pointing away from 0, due to the normal reflection of $X_t$ on the boundary of $D$. Hence, there is some $p_1 > 0$, independent of $\varepsilon < 1/200$, such that $Z_t$ may reach 1 within $1/2$ unit of time, for any starting point of $X_t$ inside $M^o$, with probability greater than $2p_1$. In other words, if $X_0 \in M^o$ then with probability greater than $2p_1$, the process $X_t$ will hit $K_3 \cup K_7$ before time $t = 1/2$. By symmetry, the process will be more likely to hit $K_3$ first, if it starts to the left of $K_5$ (i.e. in $M_i^o$), and it will be more likely to hit $K_7$ first if it starts in $M_e^o$.

The same analysis applies to the other two "bridges" of $D$. Hence, there exists $p_1 > 0$ such that for all $\varepsilon \in (0, 1/200)$, for all $x \in M_i$ and all $x' \in M_e$,

$$P(\tau_I < 1/2 \mid X_0 = x) > p_1 \quad \text{and} \quad P(\tau_E < 1/2 \mid X_0 = x') > p_1.$$

Suppose that $\gamma \subset D$ is a connected set of diameter greater than $10^{-10}$ such that $\gamma \cap I \neq \emptyset$. It is easy and elementary to prove the following: There exists $p_2 > 0$ such that for all $\varepsilon \in (0, 1/200)$, for all $x \in I$, for all $\gamma \subset D$ satisfying the above conditions,

$$P(\tau_\gamma < 1/2 \mid X_0 = x) \geq P(\tau_\gamma < 1/2, \ \tau_{\mathcal{T}K_4} > 1 \mid X_0 = x) > p_2.$$

Note that $p_2$ is independent of $\varepsilon$ as the second probability in the last formula depends only on the connected component of $D \setminus \mathcal{T}K_4$ containing $\underline{0}$ and this component is independent of $\varepsilon$.

One can also easily state and derive the counterpart of this result for the outer domain $E$.

*Step* 3. We are now ready to prove Lemma 2.

Assume first that $\Gamma \cap E = \emptyset$ and that $\Gamma \cap I \neq \emptyset$. By the Courant Nodal Line Theorem (Courant and Hilbert (1953)) the nodal line $\Gamma$ divides $D$ into two connected components. Under the current assumptions, one of these two components is a subset of $M \cup I$; we will call this component $D_c$. The first eigenvalue of the Laplacian in $D_c$ with mixed boundary conditions, Dirichlet on $\Gamma$ and Neumann elsewhere on the boundary of $D_c$ is exactly $\mu_2$.

Let $\Gamma_1$ denote a connected component of $\Gamma$ that intersects $I$. Step 1 implies that the diameter of $\Gamma_1$ is at least $10^{-10}$. Hence, using the results of Step 2, we get that for all $\varepsilon < 1/200$, for all $x \in I \cap D_c$,

$$P(\tau_\Gamma \leq 1/2 \mid X_0 = x) \geq p_2.$$

On the other hand, for all $x \in M_i$, using the strong Markov property at time $\tau_I$, Step 2 and the last inequality, we get that

$$P(\tau_\Gamma \leq 1 \mid X_0 = x) \geq p_2 p_1.$$

Finally, for all $x \in M_e \cap D_c$, as $E \cap D_c = \emptyset$,

$$P_x(\tau_\Gamma \leq 1 \mid X_0 = x) \geq P(\tau_E \leq 1 \mid X_0 = x) \geq p_1.$$

Hence, for all $x \in D_c$,

$$P(\tau_\Gamma \leq 1 \mid X_0 = x) \geq p_1 p_2,$$

so that the Markov property applied at times $n = 1, 2, \ldots$, implies that for all $n \geq 1$,

$$P(\tau_\Gamma \geq n \mid X_0 = x) \leq (1 - p_1 p_2)^n$$

and consequently that $\mu_2 \geq -\log(1 - p_1 p_2)$. Note that $p_1$ and $p_2$ are independent of $\varepsilon < 1/200$. Hence, combining this with Lemma 1 shows that for small enough $\varepsilon$, one never has $\{\Gamma \cap I \neq \emptyset$ and $\Gamma \cap E = \emptyset\}$.

The other two cases, namely $\{\Gamma \cap E \neq \emptyset$ and $\Gamma \cap I = \emptyset\}$ and $\{\Gamma \cap E \neq \emptyset$ and $\Gamma \cap I \neq \emptyset\}$, can be dealt with in the same way. Hence, for small $\varepsilon$, $\Gamma \subset M$. □

*Remark.* In almost exactly the same way, one could prove that the nodal line is in fact confined to an arbitrarily small neighbourhood of $\mathcal{T}K_6$ when $\varepsilon$ is sufficiently small, but Lemma 2 is sufficient for our purposes.

In the rest of the paper $\varepsilon > 0$ is assumed to be small enough so that the nodal line of any second Neumann eigenfunction in $D$ is a subset of $M$.

## 5. The second eigenvalue is simple

Our proof of the fact that the second eigenvalue is simple is based on an almost trivial argument. However, this argument seems to be so useful that we state it as a lemma. It originally appeared in the proofs of Propositions 2.4 and 2.5 of Bañuelos and Burdzy (1999).

LEMMA 3. *Suppose that there exists $z_0 \in D$ such that the nodal line of any second Neumann eigenfunction does not contain $z_0$. Then the second eigenvalue is simple.*

*Proof.* Suppose that $\varphi_2$ and $\widetilde{\varphi}_2$ are two independent eigenfunctions corresponding to $\mu_2$. By assumption, $\varphi_2(z_0) \neq 0$ and $\widetilde{\varphi}_2(z_0) \neq 0$ so the function

$$x \mapsto \varphi_2(x)\widetilde{\varphi}_2(z_0) - \varphi_2(z_0)\widetilde{\varphi}_2(x)$$

is a nonzero eigenfunction corresponding to $\mu_2$. Since it vanishes at $z_0$, we obtain a contradiction. □

The lemma applies to our domain $D$ because $\Gamma \subset M$.

## 6. Gradient direction for the second eigenfunction

This final part of the proof follows the arguments of Bañuelos and Burdzy (1999) so closely that we will only present a sketch and refer the reader to that paper for more details.

Let $A$ denote the disc $B(\underline{0}, 1/10)$, and let $u(t, x)$ be the solution to the Neumann heat problem in $D_1$ with the initial temperature $u(0, x) = \mathbf{1}_{D_1 \cap A}(x)$.

We set $u(t,y) = u(t,x)$ for all $y \in \mathcal{T}x$ and then extend the function $u(t,x)$ to all $x \in D$ by continuity. Due to the fact that $u$ satisfies the Neumann boundary conditions in $D_1$ and the symmetry, it is clear that $u(t,x)$ solves the Neumann heat equation in $D$ with the initial condition $u(0,x) = \mathbf{1}_A(x)$.

Since the nodal line of $\varphi_2$ is confined to $M$, the sign of $\mathbf{1}_A(x)\varphi_2(x)$ is constant. We conclude that

$$c_2 = \int_D u(0,x)\varphi_2(x)dx = \int_A \varphi_2(x)dx \neq 0,$$

and so the second eigenfunction coefficient $c_2$ is nonzero in the eigenfunction expansion for $u(t,x)$,

$$u(t,x) = c_1 + c_2\varphi_2(x)e^{-\mu_2 t} + \ldots.$$

With no loss of generality, we can assume that $c_2 > 0$, choosing the sign of $\varphi_2$ accordingly. But (see, e.g., Proposition 2.1 of Bañuelos and Burdzy (1999)),

$$u(t,x) = c_1 + c_2\varphi_2(x)e^{-\mu_2 t} + R(t,x), \qquad x \in D, \ t \geq 0,$$

where $R(t,x)$ converges to 0 as $t \to \infty$ faster than $e^{-\mu_2 t}$, uniformly in $x \in D$. Hence, if we can show that for some fixed $x, y \in D$ and all $t > 0$ we have $u(t,x) \geq u(t,y)$ then we must also have $\varphi_2(x) \geq \varphi_2(y)$.

Recall that $\alpha_1$ and $\alpha_2$ denote the minimum and maximum of the angles between vectors $\overrightarrow{A_j A_{j+1}}$, $j = 1, 2, \ldots, 9$, and the horizontal axis and that $\alpha_2 - \alpha_1 < \pi/2$. In view of this fact, the arguments of Theorems 3.1 and 3.2 (see also Example 3.2) of Bañuelos and Burdzy (1999) can be easily adjusted to our domain $D_1$ and imply that with our choice of the initial condition for $u(t,x)$, we have $u(t,x) \geq u(t,y)$ whenever the angle between the vector $\overrightarrow{xy}$ and the horizontal axis lies within $(\alpha_2 - \pi/2, \alpha_1 + \pi/2)$. Hence, we have $\varphi_2(x) \geq \varphi_2(y)$ for all such $x, y \in D_1$. In particular, for every $x \in D_1$, $\varphi_2(\underline{0}) \geq \varphi_2(x)$. In order to prove the strong inequality we observe that for every $x \in \overline{D_1} \setminus \{\underline{0}\}$, we can find a nonempty open set $F_x \subset D_1$ such that for every $y \in F_x$, the angles formed by the vectors $\overrightarrow{\underline{0}y}$ and $\overrightarrow{yx}$ with the horizontal axis belong to $(\alpha_2 - \pi/2, \alpha_1 + \pi/2)$. If $\varphi_2(\underline{0}) = \varphi_2(x)$ then $\varphi_2(\underline{0}) = \varphi_2(y) = \varphi_2(x)$ for all $y \in F_x$. The remark following Corollary (6.31) in Folland (1976) may be applied to the operator $\Delta + \mu_2$ to conclude that the eigenfunctions are real analytic and therefore they cannot be constant on an open set unless they are constant on the whole domain $D$. We see that $\varphi_2$ attains its strict maximum in $\overline{D_1}$ at the point $\underline{0}$. Since the same argument applies to every set $\sigma(D_1)$ for all $\sigma \in G$, the function $\varphi_2$ attains its strict maximum in $D$ at $\underline{0}$. This completes the proof of Theorem 1. □

University of Washington, Seattle, WA
*E-mail address*: burdzy@math.washington.edu

Université Paris-Sud, Orsay, France
*E-mail address*: wendelin.werner@math.u-psud.fr